\documentclass[a4paper,11pt]{article}

\usepackage{pdfpages}
\usepackage{graphicx}
\usepackage{mathtools}
\usepackage{verbatim}
\usepackage{amssymb}
\usepackage{bm}
\usepackage{enumitem}
\usepackage{cases}
\usepackage{latexsym}
\usepackage{amsthm}
\usepackage{amsfonts}
\usepackage{bbm}
\usepackage{etoolbox}
\usepackage[top=25mm, bottom=25mm, left=28mm, right=28mm]{geometry}
\usepackage{cite}
\usepackage{changepage}
\usepackage{fancybox}
\usepackage{tikz}

\newtheorem{theorem}{Theorem}[section]
\AfterEndEnvironment{theorem}{\noindent\ignorespaces}
\newtheorem{corollary}[theorem]{Corollary}
\AfterEndEnvironment{corollary}{\noindent\ignorespaces}
\newtheorem{lemma}[theorem]{Lemma}
\AfterEndEnvironment{lemma}{\noindent\ignorespaces}

\AfterEndEnvironment{proposition}{\noindent\ignorespaces}

\AfterEndEnvironment{question}{\noindent\ignorespaces}
\newtheorem{conjecture}[theorem]{Conjecture}
\AfterEndEnvironment{conjecture}{\noindent\ignorespaces}

\AfterEndEnvironment{problem}{\noindent\ignorespaces}

\theoremstyle{definition}
\newtheorem{definition}[theorem]{Definition}
\AfterEndEnvironment{definition}{\noindent\ignorespaces}

\theoremstyle{remark}

\AfterEndEnvironment{example}{\noindent\ignorespaces}
\newtheorem*{remark}{Remark}
\AfterEndEnvironment{remark}{\noindent\ignorespaces}

\AfterEndEnvironment{notation}{\noindent\ignorespaces}

\newenvironment{proof2}{\proof[\textnormal{\textbf{Proof.}}]}{\qed}

\usepackage{chngcntr}
\usepackage{apptools}
\AtAppendix{\counterwithin{theorem}{section}}

\def\e{\varepsilon}

\def\d{\delta}

\def\s{\subset}
\def\ex{\mathrm{ex}}

\def\cS{\mathcal{S}}
\def\cT{\mathcal{T}}

\begin{document}

\title{The extremal number of the subdivisions of the complete bipartite graph}

\author{Oliver Janzer\thanks{Department of Pure Mathematics and Mathematical Statistics, University of Cambridge, United Kingdom.
E-mail: {\tt oj224@cam.ac.uk}.}}

\date{}

\maketitle

\begin{abstract}
	
For a graph $F$, the $k$-subdivision of $F$, denoted $F^k$, is the graph obtained by replacing the edges of $F$ with internally vertex-disjoint paths of length $k$. In this paper, we prove that $\mathrm{ex}(n,K_{s,t}^k)=O(n^{1+\frac{s-1}{sk}})$, which is tight for $t$ sufficiently large. This settles a conjecture of Conlon--Janzer--Lee, and improves on a substantial body of work by Conlon--Janzer--Lee and Jiang--Qiu.

\end{abstract}

\section{Introduction}

If $H$ is a graph, we write $\ex(n,H)$ for the maximal number of edges in a graph on $n$ vertices which does not contain $H$ as a subgraph. The estimation of $\ex(n,H)$ for various graphs~$H$ is one of the main areas of research in extremal graph theory. In the case where $H$ has chromatic number at least $3$, the famous Erd\H os-Stone-Simonovits theorem \cite{ES46,ESi66} determines the asymptotics of $\ex(n,H)$, but for general bipartite graphs $H$, the function is much less understood. There are, however, several results on the extremal number of subdivided graphs.

For a graph $F$, the \emph{$k$-subdivision} $F^k$ of $F$ is the graph obtained by replacing the edges of $F$ with internally vertex-disjoint paths of length $k$. We remark that often the same graph is called the $(k-1)$-subdivision, and denoted $F^{k-1}$. The estimation of the extremal number of subdivisions has attracted the attention of many researchers recently. See \cite{Ma67,KP88,J11,JS12,CL18,Ja18,KKL18,CJL19,Ja19longer,JQ19} for results on the subject.

In this paper, we focus on the extremal number of the subdivisions of the complete bipartite graph. The first few results on this topic concerned the $2$-subdivision of the complete bipartite graph. Conlon and Lee \cite{CL18} proved that if $s\leq t$, then $\ex(n,K_{s,t}^2)=O(n^{\frac{3}{2}-\frac{1}{12t}})$. (Here, and everywhere else in the paper, it is assumed that $n\rightarrow \infty$ and other parameters are kept constant. In particular, the implied constants may depend on $s$, $t$ and $k$.) This was improved by the author \cite{Ja18} to $\ex(n,K_{s,t}^2)=O(n^{\frac{3}{2}-\frac{1}{4s-2}})$ and the same result was reproved using different methods by Kang, Kim and Liu \cite{KKL18}. Moreover, they conjectured that $\ex(n,K_{s,t}^2)=O(n^{\frac{3}{2}-\frac{1}{2s}})$ holds, which is tight for sufficiently large $t$ by a general result of Bukh and Conlon~\cite{BC17} (see Theorem \ref{thmBC} below). The conjecture was proved by Conlon, Janzer and Lee \cite{CJL19}. About longer subdivisions, they proved the following result.

\begin{theorem}[Conlon--Janzer--Lee \cite{CJL19}] \label{bipsubnotsharp}
	For any integers $s,t,k\geq 2$, $$\ex(n,K_{s,t}^k)=O(n^{1+\frac{s}{sk+1}}).$$
\end{theorem}

This is nearly sharp for $t$ sufficiently large, since the above result of Bukh and Conlon \cite{BC17} implies that there exists $t_0=t_0(s,k)$ such that for all $t\geq t_0$, $\ex(n,K_{s,t}^k)=\Omega(n^{1+\frac{s-1}{sk}})$.

Conlon, Janzer and Lee conjectured that this lower bound is tight.

\begin{conjecture}[Conlon--Janzer--Lee \cite{CJL19}] \label{conjecturebipsub}
	For any integers $s,t,k\geq 2$, $$\ex(n,K_{s,t}^k)=O(n^{1+\frac{s-1}{sk}}).$$
\end{conjecture}

Jiang and Qiu proved that the conjecture holds for $k=3$ and $k=4$.

\begin{theorem}[Jiang--Qiu \cite{JQ19}]
	For any integers $s,t\geq 2$ and $k\in \{3,4\}$, $$\ex(n,K_{s,t}^k)=O(n^{1+\frac{s-1}{sk}}).$$
\end{theorem}

In this paper we prove Conjecture \ref{conjecturebipsub} for arbitrary $k$.

\begin{theorem} \label{bipsub}
	For any integers $s,t,k\geq 2$, $$\ex(n,K_{s,t}^k)=O(n^{1+\frac{s-1}{sk}}).$$
\end{theorem}

As mentioned above, this is tight for $t$ sufficiently large.

\begin{corollary} \label{exactbip}
	For any integers $s,k\geq 2$, there exists $t_0=t_0(s,k)$ such that for all integers $t\geq t_0$, $$\ex(n,K_{s,t}^k)=\Theta(n^{1+\frac{s-1}{sk}}).$$
\end{corollary}

In the next section we state our main lemma and deduce Theorem \ref{bipsub} from it. In Section~\ref{sec3bip} we prove the main lemma. Section \ref{sec4bip} contains some concluding remarks.

\section{The proof of Theorem \ref{bipsub} given the main lemma}

Let us start with a standard lemma that allows us to restrict our attention to nearly regular host graphs. The first variant of this lemma was due to Erd\H os and Simonovits \cite{ES70}, but we will use a version given by Jiang and Seiver \cite{JS12}.

A graph $G$ is called $K$-almost-regular if $\max_{v\in V(G)} d(v)\leq K\min_{v\in V(G)} d(v)$, where $d(v)$ is the degree of vertex $v$.

\begin{lemma}[Jiang--Seiver \cite{JS12}] \label{lemmaJS}
	Let $\e,c$ be positive reals, where $\e<1$ and $c\geq 1$. Let $n$ be a positive integer that is sufficiently large as a function of $\e$. Let $G$ be a graph on $n$ vertices with $e(G)\geq cn^{1+\e}$. Then $G$ contains a $K$-almost-regular subgraph $G'$ on $m\geq n^{\frac{\e-\e^2}{2+2\e}}$ vertices such that $e(G')\geq \frac{2c}{5}m^{1+\e}$ and $K=20\cdot 2^{\frac{1}{\e^2}+1}$.
\end{lemma}

Using this lemma, Theorem \ref{bipsub} reduces to the following statement.

\begin{theorem} \label{bipsubreduced}
	Let $s,t,k\geq 2$ be integers. Let $G$ be a $K$-almost-regular graph on $n$ vertices with minimum degree $\d=\omega(n^{\frac{s-1}{sk}})$. Then, for $n$ sufficiently large, $G$ contains $K_{s,t}^k$ as a subgraph.
\end{theorem}

In what follows, let us fix the integers $s,t,k\geq 2$. It will be tacitly assumed throughout the paper that $n$ is sufficiently large compared to all other parameters.

The next definition is due to Jiang and Qiu \cite{JQ19}.

\begin{definition}
	Let $\ell_1,\dots,\ell_s$ be positive integers. An \emph{$s$-legged spider} $S$ with length vector $(\ell_1,\dots,\ell_s)$ consists of a vertex $u$, called the centre of the spider, and paths $P_1,\dots,P_s$, called the legs of $S$, of lengths $\ell_1,\dots,\ell_s$, starting at $u$ and sharing no vertex other than $u$. For convenience, we define two spiders $S$ and $S'$ to be different if $P_i\neq P'_i$ for some $1\leq i\leq s$, where $P'_1,\dots,P'_s$ are the legs of $S'$. So different spiders can form the same graph, e.g. if $\ell_1=\ell_2$, $P_1=P'_2$, $P_2=P'_1$ and $P_i=P'_i$ for $i\geq 3$.
	
	Let $v_i$ be the endpoint of $P_i$ different from $u$. Then we say that $S$ has leaf vector $(v_1,\dots,v_s)$.	
	
	We say that $S'$ is a subspider of $S$ if they have the same centre and for each $1\leq i\leq s$, the $i$th leg of $S'$ is a subpath of the $i$th leg of $S$.
\end{definition}

To prove Theorem \ref{bipsubnotsharp}, Conlon, Janzer and Lee showed in \cite{CJL19} that (roughly speaking) if a graph has many pairs of short paths $(P,P')$ such that $P$ and $P'$ are of equal length and have the same endpoints, then the graph contains $K_{s,t}^k$ as a subgraph. In this paper we shall prove an analogous statement for spiders; that is, if there are many pairs of spiders $(S,S')$ such that $S$ and $S'$ have the same length vector and the same leaf vector, then the graph contains $K_{s,t}^k$ as a subgraph. We remark that Jiang and Qiu used a similar strategy in \cite{JQ19}. The main contribution of this paper is that we have a different approach for finding copies of $K_{s,t}^k$ given a collection of spiders with a fixed length vector. This allows us to deal with spiders with arbitrary length vectors (see Lemma \ref{notgoodspiders} below), whereas Jiang and Qiu only deal with a certain subset of all possible length vectors which suffices in the cases $k=3,4$.

The next definition extends Definition 6.2 from \cite{CJL19} to spiders.

\begin{definition}
	Let $L\geq 1$ be a real number. Let $f(1,L)=L$, and for $2\leq \ell\leq sk$, let $$f(\ell,L)=1+f(\ell-1,L)^{16}(\ell-1)^2\max_{1\leq i\leq \ell-1}f(i,L)f(\ell-i,L).$$
	
	We define the notions of \emph{$L$-admissible} and \emph{$L$-good} paths and spiders recursively as follows.
	
	Every path of length $1$ is both $L$-admissible and $L$-good. For $2\leq \ell\leq k$, we say that a path of length $\ell$ is $L$-admissible if each of its proper subpaths is $L$-good. A path of length~$\ell$ is $L$-good if it is $L$-admissible and the number of $L$-admissible paths of length $\ell$ between their endpoints is at most $f(\ell,L)$.
	
	Every $s$-legged spider with length vector $(1,\dots,1)$ is $L$-admissible. Now let $1\leq \ell_1,\dots,\ell_s\leq k$ and assume that $\ell_i>1$ for some $i$. A spider with centre $u$ and legs $P_i=uw_{i,1}\dots w_{i,\ell_i}$ (for $1\leq i\leq s$) is $L$-admissible if the following two conditions hold:
	\begin{itemize}
		\item for any $1\leq i\leq s$ and any $1\leq j<\ell_i$, the $s$-legged spider with centre $u$ and legs~$P_1,\dots,P_{i-1},P'_i,P_{i+1},\dots,P_s$ is $L$-good, where $P'_i=uw_{i,1}\dots w_{i,j}$
		
		\item for any $1\leq i\leq s$, the path $P_i$ is $L$-good.
	\end{itemize}
	Finally, we say that a spider with length vector $(\ell_1,\dots,\ell_s)$ and leaf vector $(v_1,\dots,v_s)$ is $L$-good if it is $L$-admissible and the number of $L$-admissible spiders with length vector $(\ell_1,\dots,\ell_s)$ and leaf vector $(v_1,\dots,v_s)$ is at most $f(\ell,L)$, where $\ell=\ell_1+\dots+\ell_s$.
\end{definition}

\begin{remark}
	\begin{enumerate}[label=(\arabic*)]
		\item This is well-defined since whether a spider is $L$-admissible or not depends only on the $L$-goodness of smaller spiders and paths.
		\item $L$ will be chosen to be a constant not depending on $n$.
	\end{enumerate}	
\end{remark}

\noindent
As we mentioned earlier, paths have already been satisfyingly controlled in \cite{CJL19}. We will use the next lemma, which follows easily from Corollary 6.9 in \cite{CJL19} since the graph $H$ there contains $K_{s,t}^k$ as a subgraph.

\begin{lemma}[Conlon--Janzer--Lee \cite{CJL19}] \label{kpathsgood}
	Let $G$ be a $K_{s,t}^{k}$-free $K$-almost-regular graph on $n$ vertices with minimum degree~$\d=\omega(1)$. Then for any $1\leq j\leq k$, the number of paths of length $j$ which are not $L$-good is at most $c_Ln\d^j$, where $c_L\rightarrow 0$ as $L\rightarrow \infty$.
\end{lemma}

The main technical result of this paper is the following lemma, which is the analogue of Lemma \ref{kpathsgood} for spiders.

\begin{lemma} \label{notgoodspiders}
	Let $G$ be a $K_{s,t}^{k}$-free $K$-almost-regular graph on $n$ vertices with minimum degree~$\d=\omega(1)$ and let $1\leq \ell_1,\dots,\ell_s\leq k$. Then the number of $s$-legged spiders with length vector $(\ell_1,\dots,\ell_s)$ which are $L$-admissible but not $L$-good is at most $c'_Ln\d^{\ell_1+\dots +\ell_s}$, where $c'_L\rightarrow 0$ as $L\rightarrow \infty$.
\end{lemma}

We postpone the proof of this lemma to the next section and first show how it implies Theorem \ref{bipsubreduced}. The next lemma is an easy corollary of Lemma \ref{notgoodspiders}.

\begin{lemma} \label{kspidersgood}
	Let $G$ be a $K_{s,t}^{k}$-free $K$-almost-regular graph on $n$ vertices with minimum degree~$\d=\omega(1)$. Then the number of $s$-legged spiders with length vector $(k,k\dots,k)$ which are not $L$-good is at most $c''_Ln\d^{sk}$, where $c''_L\rightarrow 0$ as $L\rightarrow \infty$.
\end{lemma}

\begin{proof2}
	Suppose that some $s$-legged spider $S$ with length vector $(k,\dots,k)$ and legs $P_1,\dots,P_s$ is not $L$-good.
	
	We distinguish two cases. First, assume that some $P_i$ is not $L$-good. By Lemma \ref{kpathsgood}, there are at most $c_Ln\d^k$ choices for $P_i$, where $c_L\rightarrow 0$ as $L\rightarrow \infty$. Since the maximum degree of $G$ is at most $K\d$, the number of ways to extend a given $P_i$ to an $s$-legged spider with length vector $(k,\dots,k)$ is at most $(K\d)^{(s-1)k}$. Thus, the number of $s$-legged spiders with length vector $(k,\dots,k)$ such that one of the legs is not $L$-good is at most $s\cdot c_Ln\d^k\cdot (K\d)^{(s-1)k}=sK^{(s-1)k}c_Ln\d^{sk}$.
	
	Now assume that all the $P_i$ are $L$-good. Choose an $s$-legged subspider $S'$ with the same centre and legs $P'_1,\dots,P'_s$ which are subpaths of $P_1,\dots,P_s$ such that $S'$ is minimal with respect to the condition that $S'$ is not $L$-good. Let $\ell_i$ be the length of $P'_i$. Suppose that $S'$ is not $L$-admissible. Since each $P_i$ is $L$-good, so is every subpath of every leg $P'_i$. Thus, there must be a proper $s$-legged subspider in $S'$ which is not $L$-good. This contradicts the minimality of $S'$. So $S'$ is $L$-admissible but not $L$-good. By Lemma \ref{notgoodspiders}, for any fixed $1\leq \ell_1,\dots,\ell_s\leq k$, the number of $s$-legged spiders with length vector $(k,\dots,k)$ whose subspider with length vector $(\ell_1,\dots,\ell_s)$ is $L$-admissible but not $L$-good is at most $c'_Ln\d^{\ell_1+\dots+\ell_s}\cdot(K\d)^{sk-\ell_1-\dots-\ell_s}$. Summing over all choices for $\ell_1,\dots,\ell_s$, we find that the number of $s$-legged spiders with length vector $(k,\dots,k)$ which are not $L$-good but whose legs are all $L$-good is at most $k^s\cdot K^{sk} c'_Ln\d^{sk}$.
\end{proof2}

\medskip

We are now in a position to complete the proof of Theorem \ref{bipsubreduced}.

\begin{proof}[\textnormal{\textbf{Proof of Theorem \ref{bipsubreduced}}}]
	Choose $L$ such that the $c''_L$ provided by Lemma \ref{kspidersgood} satisfies $c''_L\leq 1/2$. Then by Lemma \ref{kspidersgood}, for $n$ sufficiently large, the number of $L$-good $s$-legged spiders with length vector $(k,\dots,k)$ is at least $\frac{1}{3}n\d^{sk}>f(sk,L)n^s$. Thus, there exists an $s$-tuple $(v_1,\dots,v_s)$ of vertices such that the number of $L$-good $s$-legged spiders with length vector $(k,\dots,k)$ and leaf vector $(v_1,\dots,v_s)$ is greater than $f(sk,L)$. This contradicts the definition of an $L$-good spider.
\end{proof}

\section{The proof of Lemma \ref{notgoodspiders}} \label{sec3bip}

In this section we prove Lemma \ref{notgoodspiders}, after which the proof of our main theorem is complete. For this section, we fix some $1\leq \ell_1,\dots,\ell_s\leq k$ and write $\ell=\ell_1+\dots+\ell_s$.

In what follows, it will be crucial to look at "spiders" some of whose legs may consist of zero edges.

\begin{definition}
	Let $\ell'_1,\dots,\ell'_s$ be nonnegative integers. A \emph{generalised spider} $S$ with length vector $(\ell'_1,\dots,\ell'_s)$ consists of a vertex $u$ (the centre of $S$) and paths $P_1,\dots,P_s$ (the legs of $S$) of lengths $\ell'_1,\dots,\ell'_s$, starting at $u$ and sharing no vertex other than $u$. Let $P_i$ have endpoints $u$ and $v_i$. Then we say that $S$ has leaf vector $(v_1,\dots,v_s)$.
\end{definition}

The next lemma states that if there are many $L$-admissible but not $L$-good spiders with length vector $(\ell_1,\dots,\ell_s)$ in our graph, then we can find many $L$-admissible spiders with length vector $(\ell_1,\dots,\ell_s)$ and some useful extra properties.

\begin{lemma} \label{regularset}
	Let $G$ be a $K$-almost-regular graph on $n$ vertices with minimum degree $\d$. Assume that $L$ is sufficiently large compared to $s$, $k$ and $K$ and that there are at least $\frac{n\d^{\ell_1+\dots+\ell_s}}{L}$ $L$-admissible but not $L$-good spiders with length vector $(\ell_1,\dots,\ell_s)$. Then there exists a non-empty set $\cS$ of $L$-admissible spiders with length vector $(\ell_1,\dots,\ell_s)$ such that the following conditions hold.
	\begin{enumerate}[label=(\roman*)]
		\item For any $S\in \cS$, the number of spiders $T\in \cS$ with the same leaf vector as that of $S$ is at least $\frac{f(\ell,L)}{2}$.
		\item For any $S\in \cS$, and any $\gamma_1,\dots,\gamma_s\in \{0,1\}$, the subspider of $S$ with length vector $(\ell_1-\gamma_1,\dots,\ell_s-\gamma_s)$ (which is a generalised spider) is contained as a subspider in at least $\frac{\d^{\gamma_1+\dots+\gamma_s}}{L^2}$ elements of $\cS$. 
	\end{enumerate}
\end{lemma}

\begin{proof2}
	Define a sequence of sets $\cT_0,\cT_1,\dots,\cT_m$ recursively as follows. Take $\cT_0$ be the set of all $L$-admissible but not $L$-good spiders with length vector $(\ell_1,\dots,\ell_s)$. Then, if there is some $S\in \cT_i$ which violates condition (i), ie. the number of spiders $T\in \cT_i$ with the same leaf vector as that of $S$ is less than $\frac{f(\ell,L)}{2}$, then choose such an $S$ arbitrarily and let $\cT_{i+1}=\cT_i\setminus \{S\}$. Also, if no such $S$ exists, but there is some $S\in \cT_i$ which violates condition (ii), ie. there exist some $\gamma_1,\dots,\gamma_s\in \{0,1\}$ such that the subspider of $S$ with length vector $(\ell_1-\gamma_1,\dots,\ell_s-\gamma_s)$ is contained in less than $\frac{\d^{\gamma_1+\dots+\gamma_s}}{L^2}$ elements of $\cT_i$, then choose such an $S$ arbitrarily and let $\cT_{i+1}=\cT_i\setminus \{S\}$. The process eventually terminates with some set $\cT_m$. Let $\cS=\cT_m$. It is clear that $\cS$ satisfies conditions (i) and (ii); all we need to prove is that $\cS\neq \emptyset$. Note that every $S\in \cT_0$ is $L$-admissible but not $L$-good, so there are at least $f(\ell,L)$ elements $T\in \cT_0$ with the same leaf vector as that of $S$. Among the set of elements of $\cT_0$ with a fixed leaf vector, at most $\frac{f(\ell,L)}{2}$ are discarded because of violating condition (i) at some point. Thus, if $\cS=\emptyset$, then at least half of the elements of $\cT_0$, and so at least $\frac{n\d^{\ell_1+\dots+\ell_s}}{2L}$ spiders are discarded because of violating condition (ii) at some point. However, any generalised spider $R$ with length vector $(\ell_1-\gamma_1,\dots,\ell_s-\gamma_s)$ is "responsible" for discarding at most $\frac{\d^{\gamma_1+\dots+\gamma_s}}{L^2}$ elements, meaning that the number of elements discarded because they contain $R$ which is contained in less than $\frac{\d^{\gamma_1+\dots+\gamma_s}}{L^2}$ elements of some $\cT_i$ is at most $\frac{\d^{\gamma_1+\dots+\gamma_s}}{L^2}$. Since the number of generalised spiders with length vector $(\ell_1-\gamma_1,\dots,\ell_s-\gamma_s)$ is at most $n(K\d)^{(\ell_1-\gamma_1)+\dots+(\ell_s-\gamma_s)}$, the total number of elements discarded because of violating condition (ii) at some point is at most $2^s\cdot \frac{n(K\d)^{\ell_1+\dots+\ell_s}}{L^2}$. For $L>2^{s+1} K^{\ell_1+\dots+\ell_s}$, this is less than $\frac{n\d^{\ell_1+\dots+\ell_s}}{2L}$, contradicting our earlier claim. Thus, $\cS\neq \emptyset$.
\end{proof2}

\begin{lemma} \label{disjointspiders}
	Let $L\geq 1$ be real and let $v_1,\dots,v_s$ be vertices.
	Suppose that there is a set $\cT$ of at least $\frac{f(\ell,L)}{2}$ $L$-admissible spiders with length vector $(\ell_1,\dots,\ell_s)$ and leaf vector $(v_1,\dots,v_s)$. Then, among these, there exist at least $\frac{f(\ell-1,L)^{16}}{2}$ spiders which are pairwise vertex-disjoint apart from their leaves.
\end{lemma}

\begin{proof2}
	Suppose otherwise. Take a maximal set of such spiders. By assumption, we have chosen at most $\frac{f(\ell-1,L)^{16}}{2}$ spiders. Each such spider has $\ell+1-s\leq \ell-1$ non-leaf vertices, so altogether they have at most $\frac{f(\ell-1,L)^{16}(\ell-1)}{2}$ non-leaf vertices. By the maximality assumption, each $S\in \cT$ contains at least one of these vertices. Thus, by the pigeonhole principle, there exist some vertex $x$ and a set $\mathcal{S}\s \cT$ of size at least $\frac{f(\ell,L)/2}{(\ell-1)\cdot f(\ell-1,L)^{16}(\ell-1)/2}$ such that the elements of $\cS$ all contain the vertex $x$ in the same non-leaf position (meaning that there are $i$ and $j<\ell_i$ such that in all $S\in \mathcal{S}$, $x$ is the $j$th vertex on the $i$th leg, where the centre of the spider is viewed as the $0$th vertex on the leg). Note that $|\mathcal{S}|\geq  \frac{f(\ell,L)}{f(\ell-1,L)^{16}(\ell-1)^2}>\max_{1\leq b\leq \ell-1} f(b,L)f(\ell-b,L)$.
	
	We now distinguish two cases. First, let us assume that $x$ is not the centre in the spiders in $\mathcal{S}$. Then there exists some $1\leq i\leq s$ and some $1\leq j< \ell_i$ such that $x$ is the $j$th vertex on the $i$th leg in each of these spiders. Let $b=\ell_i-j$. Since $|\cS|>f(b,L)f(\ell-b,L)$ and each element of $\cS$ is $L$-admissible, either there are more than $f(b,L)$ $L$-good paths of length $b$ between $x$ and $v_i$ or there are more than $f(\ell-b,L)$ $L$-good $s$-legged spiders with length vector $(\ell_1,\dots,\ell_{i-1},j,\ell_{i+1},\dots,\ell_s)$ and leaf vector $(v_1,\dots,v_{i-1},x,v_{i+1},\dots,v_s)$. The first case contradicts the definition of an $L$-good path and the second case contradicts the definition of an $L$-good $s$-legged spider.
	
	Let us now assume that $x$ is the centre in the spiders in $\mathcal{S}$. Note that $$|\mathcal{S}|>\max_{1\leq b\leq \ell-1} f(b,L)f(\ell-b,L)\geq f(\ell_1,L)f(\ell-\ell_1,L)\geq f(\ell_1,L)f(\ell_2,L)\dots f(\ell_s,L),$$
	where the last inequality follows easily from the recursive definition of $f$. Thus, there exists some $i\leq s$ such that there are more than $f(\ell_i,L)$ $L$-good paths of length $\ell_i$ between $x$ and $v_i$. This contradicts the definition of an $L$-good path.
\end{proof2}

\medskip

In the key part of the proof of Lemma \ref{notgoodspiders} it will be necessary to assume that $\ell_i=1$ holds for at most one choice of $i$. Accordingly, we first deal with the other case separately.

\begin{lemma} \label{length11}
	Let $G$ be a $K_{s,t}^{k}$-free $K$-almost-regular graph on $n$ vertices with minimum degree~$\d=\omega(1)$, and assume that $\ell_1=\ell_2=1$. Then the number of $s$-legged spiders with length vector $(\ell_1,\dots,\ell_s)$ which are $L$-admissible but not $L$-good is at most $c'_Ln\d^{\ell_1+\dots +\ell_s}$, where $c'_L\rightarrow 0$ as $L\rightarrow \infty$.
\end{lemma}

\begin{proof2}
	If $s=2$, then the result follows from Lemma \ref{kpathsgood}, since a spider with length vector $(1,1)$ is $L$-good if and only if it is $L$-good when viewed as a path of length $2$. Assume that $s\geq 3$. Note that in this case $\ell\geq 3$.
	
	Let $S$ be an $L$-admissible but not $L$-good spider with length vector $(\ell_1,\dots,\ell_s)$ and leaf vector $(v_1,\dots,v_s)$. By definition, there exist at least $f(\ell,L)$ $L$-admissible spiders with length vector $(\ell_1,\dots,\ell_s)$ and leaf vector $(v_1,\dots,v_s)$. Hence, by Lemma \ref{disjointspiders}, for $L$ sufficiently large there exist more than $f(\ell-1,L)$ $L$-admissible spiders with length vector $(\ell_1,\dots,\ell_s)$ and leaf vector $(v_1,\dots,v_s)$ which are pairwise vertex-disjoint apart from at their leaves. In particular, there are more than $f(\ell-1,L)\geq f(2,L)$ paths of length $2$ between $v_1$ and $v_2$. Note that any path of length $2$ is $L$-admissible. Let $u$ be the centre of $S$. Then the path $v_1uv_2$ is not $L$-good.
	
	The number of ways to extend a path $xyz$ to a spider with length vector $(\ell_1,\dots,\ell_s)$, centre $y$ and first two legs $yx$ and $yz$ in this order is at most $(K\d)^{\ell_3+\dots+\ell_s}$. Thus, by Lemma~\ref{kpathsgood}, the number of $L$-admissible but not $L$-good spiders with length vector $(\ell_1,\dots,\ell_s)$ is at most $c_Ln\d^2\cdot 2 \cdot (K\d)^{\ell_3+\dots+\ell_s}$ with $c_L\rightarrow 0$ as $L\rightarrow \infty$, where the factor $c_Ln\d^2$ bounds the number of not $L$-good paths of length $2$, the factor $2$ accounts for the two edges in this path that we can use as the first leg of the spider, and the factor $(K\d)^{\ell_3+\dots+\ell_s}$ bounds the number of ways to get a spider with fixed first two legs. Since $c_Ln\d^2\cdot 2 \cdot (K\d)^{\ell_3+\dots+\ell_s}=2K^{\ell_3+\dots+\ell_s}c_Ln\d^{\ell_1+\dots+\ell_s}$, the result follows.
\end{proof2}

\medskip

Using Lemma \ref{length11} and symmetry, it is enough to prove Lemma \ref{notgoodspiders} in the case where $\ell_i=1$ holds for at most one value of $i$.

The next result is the key step in the proof of Lemma \ref{notgoodspiders}, and contains the main idea of this paper. It is proved in greater generality than what is needed for our main result, to allow for use in future work.

\begin{lemma} \label{pathsbyspiders}
	Let $\ell_i\leq k_i\leq k$ for each $i$. Assume that $\ell_i=1$ holds for at most one value of $i$. Let $G$ be a $K$-almost-regular graph on $n$ vertices with minimum degree $\d=\omega(1)$. Assume that $L$ is sufficiently large compared to $s$, $k$ and $K$ and that there exists a set $\cS$ of spiders satisfying the conditions in Lemma \ref{regularset}. For each $1\leq i\leq s$, let $\gamma_{i,0}\in \{0,1\}$ such that $k_i-\ell_i-\gamma_{i,0}$ is even. Let $R_0$ be the subspider with length vector $(\ell_1-\gamma_{1,0},\dots,\ell_s-\gamma_{s,0})$ of an arbitrary element of $\cS$. Let $R_0$ have leaf vector $(v_1,\dots,v_s)$. Let $Z\s V(G)$ be a set of size at most $L$, disjoint from $\{v_1,\dots,v_s\}$. Then there exist vertices $w_1,\dots,w_s$ and paths $P_1,\dots,P_s$ such that
	\begin{enumerate}[label=(\arabic*)]
		\item for each $i$, $P_i$ is a path of length $k_i-\ell_i$ between $v_i$ and $w_i$
		\item $(w_1,\dots,w_s)$ is the leaf vector of an element of $\cS$ and
		\item the paths $P_1,\dots,P_s$ are pairwise vertex-disjoint and avoid $Z$.	
	\end{enumerate}
\end{lemma}

\begin{proof2}
	Since $k_i-\ell_i-\gamma_{i,0}$ is an even number between $0$ and $k$, there exist $\gamma_{i,1},\dots,\gamma_{i,k-1}\in \{0,1\}$ such that $k_i-\ell_i-\gamma_{i,0}=2\gamma_{i,1}+\dots+2\gamma_{i,k-1}$.
	
	We now define a sequence $R_1,\dots,R_{k-1}$ of generalised spiders, and sequences $S_1,\dots,S_k$ and $T_1,\dots,T_k$ of spiders recursively.
	
	$R_0$ is given as a subspider of some element of $\cS$, so by property (ii) in Lemma \ref{regularset}, the number of elements of $\cS$ containing $R_0$ as a subspider is at least $\frac{\d^{\gamma_{1,0}+\dots+\gamma_{s,0}}}{L^2}$. Thus, there is some $S_1\in \cS$ containing $R_0$ such that $V(S_1)\setminus V(R_0)$ is disjoint from $Z$. Indeed, any fixed vertex not in $V(R_0)$ is a vertex in $O(\d^{\gamma_{1,0}+\dots+\gamma_{s,0}-1})$ elements of $\cS$ containing $R_0$, so the number of elements of $\cS$ containing $R_0$ and intersecting $Z\setminus V(R_0)$ is $O(\d^{\gamma_{1,0}+\dots+\gamma_{s,0}-1})$. Hence, as $\d=\omega(1)$ and $L=O(1)$, a suitable $S_1\in \cS$ indeed exists.
	
	Now choose $T_1\in \cS$ with the same leaf vector as that of $S_1$ such that $T_1$ and $S_1$ are disjoint apart from their leaves. This is possible, if $L$ is sufficiently large, by property (i) in Lemma \ref{regularset} and Lemma \ref{disjointspiders}. Let $R_1$ be the subspider of $T_1$ with length vector $(\ell_1-\gamma_{1,1},\dots,\ell_s-\gamma_{s,1})$.
	
	More generally, for any $1\leq j\leq k$, given a generalised spider $R_{j-1}$ with length vector $(\ell_1-\gamma_{1,j-1},\dots,\ell_s-\gamma_{s,j-1})$ which is a subspider of an element of $\cS$, we define $S_j$, $T_j$ and $R_j$ as follows.
	
	Choose some $S_{j}\in \cS$ containing $R_{j-1}$ such that $V(S_j)\setminus V(R_{j-1})$ is disjoint from $Z \cup (V(S_1)\cup \dots \cup V(S_{j-1})) \cup (V(T_1)\cup \dots \cup V(T_{j-1}))$. This is possible by property (ii) in Lemma~\ref{regularset}.
	
	Also, choose $T_j\in \cS$ with the same leaf vector as that of $S_j$ such that $T_j$ is disjoint from $Z \cup (V(S_1)\cup \dots \cup V(S_{j})) \cup (V(T_1)\cup \dots \cup V(T_{j-1}))$ apart from its leaves. This is possible by property (i) in Lemma \ref{regularset} and Lemma \ref{disjointspiders}.
	
	Finally, if $j<k$, let $R_j$ be the subspider of $T_j$ with length vector $(\ell_1-\gamma_{1,j},\dots,\ell_s-\gamma_{s,j})$.
	
	Now for $1\leq i\leq s$ and $0\leq j\leq k-1$, let $x_{i,2j}$ be the endpoint of the $i$th leg of $R_j$ and let $x_{i,2j+1}$ be the endpoint of the $i$th leg of $S_{j+1}$. Then, when we ignore the repetitions, the vertices $x_{i,0},x_{i,1},\dots,x_{i,2k-1}$ form a path of length $\gamma_{i,0}+2\gamma_{i,1}+\dots+2\gamma_{i,k-1}=k_i-\ell_i$. Indeed, if $\gamma_{i,0}=0$, then $x_{i,1}=x_{i,0}$ and if $\gamma_{i,0}=1$, then $x_{i,1}$ is a neighbour of $x_{i,0}$. Moreover, for any $1\leq j\leq k-1$, if $\gamma_{i,j}=0$, then $x_{i,2j+1}=x_{i,2j}=x_{i,2j-1}$ and if $\gamma_{i,j}=1$, then $x_{i,2j}$ is a neighbour of $x_{i,2j-1}$ and does not belong to $\{x_{p,q}: 1\leq p\leq s, 0\leq q\leq 2j-1\}\cup Z$, and $x_{i,2j+1}$ is a neighbour of $x_{i,2j}$ and does not belong to $\{x_{p,q}: 1\leq p\leq s, 0\leq q\leq 2j\}\cup Z$. Let $P_i$ be the path formed by the vertices $x_{i,0},x_{i,1}\dots,x_{i,2k-1}$ and let $w_i=x_{i,2k-1}$.
	
	Note that $(x_{1,0},\dots,x_{s,0})$ is the leaf vector of $R_0$, so $x_{i,0}=v_i$, therefore condition (1) in this lemma is satisfied. Moreover, $(w_1,\dots,w_s)=(x_{1,2k-1},\dots,x_{s,2k-1})$ is the leaf vector of $S_{k}$, so property (2) is also satisfied.
	
	By assumption, $Z$ is disjoint from $\{v_1,\dots,v_s\}=\{x_{1,0},\dots,x_{s,0}\}$, so it follows from the above that $P_1,\dots,P_s$ avoid $Z$. Finally, it is clear by the above discussion that if $P_1,\dots,P_s$ are not pairwise vertex-disjoint, then $x_{i,j}=x_{i',j}$ holds for some $i\neq i'$ and some $0\leq j\leq 2k-1$. However, for each $0\leq j\leq 2k-1$, $(x_{1,j},\dots,x_{s,j})$ is the leaf vector of a generalised spider whose $i$th leg consists of at least $\ell_i-1$ edges, so at most one of its legs has $0$ edges. Thus, the vertices $x_{1,j},\dots,x_{s,j}$ are distinct and condition (3) is satisfied.	
\end{proof2}

\medskip

It is not hard to connect the paths given by the previous lemma to form spiders with length vector $(k_1,\dots,k_s)$.

\begin{lemma} \label{connectthepaths}
	Let $\ell_i\leq k_i\leq k$ for each $i$. Assume that $\ell_i=1$ holds for at most one value of $i$. Let $G$ be a $K$-almost-regular graph on $n$ vertices with minimum degree $\d=\omega(1)$. Assume that $L$ is sufficiently large compared to $s$, $k$ and $K$ and that there exists a set $\cS$ of spiders satisfying the conditions in Lemma \ref{regularset}. For each $1\leq i\leq s$, let $\gamma_{i,0}\in \{0,1\}$ such that $k_i-\ell_i-\gamma_{i,0}$ is even. Let $R_0$ be the subspider with length vector $(\ell_1-\gamma_{1,0},\dots,\ell_s-\gamma_{s,0})$ of an arbitrary element of $\cS$. Let $R_0$ have leaf vector $(v_1,\dots,v_s)$. Let $Z\s V(G)$ be a set of size at most $L$, disjoint from $\{v_1,\dots,v_s\}$.
	
	Then there exists an $s$-legged spider with length vector $(k_1,\dots,k_s)$ and leaf vector $(v_1,\dots,v_s)$ that avoids $Z$.
\end{lemma}

\begin{proof2}
	Choose vertices $w_1,\dots,w_s$ and paths $P_1,\dots,P_s$ as in the conclusion of Lemma~\ref{pathsbyspiders}. $(w_{1},\dots,w_{s})$ is the leaf vector of an element of $\cS$, so by condition (i) in Lemma \ref{regularset} and Lemma \ref{disjointspiders}, there exist at least $f(\ell-1,L)^{16}/2$ spiders with length vector $(\ell_1,\dots,\ell_s)$ and leaf vector $(w_{1},\dots,w_{s})$ which are pairwise vertex-disjoint apart from at their leaves. Thus, if $L$ is sufficiently large, then there exists a spider $S$ with length vector $(\ell_1,\dots,\ell_s)$ and leaf vector $(w_{1},\dots,w_{s})$ such that $V(S)$ is disjoint from $Z$ and intersects $\bigcup_{1\leq i\leq s} V(P_{i})$ only at $\{w_{1},\dots,w_{s}\}$. Let $J_i$ be the $i$th leg of $S$, let $u$ be the centre of $S$ and let $Q_i$ be the union of $J_i$ and $P_i$. Then the spider with centre $u$ and legs $Q_1,\dots,Q_s$ is suitable.
\end{proof2}

\medskip

The next result, together with Lemma \ref{length11}, completes the proof of Lemma \ref{notgoodspiders}.

\begin{lemma} \label{findkstk}
	Assume that $\ell_i=1$ holds for at most one value of $i$. Let $G$ be a $K$-almost-regular graph on $n$ vertices with minimum degree $\d=\omega(1)$. Assume that $L$ is sufficiently large compared to $s$, $t$, $k$ and $K$ and that there are at least $\frac{n\d^{\ell_1+\dots+\ell_s}}{L}$ $L$-admissible but not $L$-good spiders with length vector $(\ell_1,\dots,\ell_s)$. Then $G$ contains $K_{s,t}^k$ as a subgraph.
\end{lemma}

\begin{proof2}
	Choose a set $\cS$ with the properties described in Lemma \ref{regularset}. Take $k_1=\dots=k_s=k$ and define $v_1,\dots,v_s$ as in the statement of Lemma \ref{connectthepaths}.
	We may repeatedly apply Lemma \ref{connectthepaths} to find $s$-legged spiders $S_1,\dots,S_t$, each with length vector $(k,\dots,k)$ and leaf vector $(v_1,\dots,v_s)$ such that $V(S_j)$ is disjoint from $(\bigcup_{1\leq i\leq j-1} V(S_i))\setminus \{v_1,\dots,v_s\}$. Then the union of these spiders is a copy of $K_{s,t}^k$.
\end{proof2}

\section{Concluding remarks} \label{sec4bip}

Let $F$ be a graph with a set $R\subsetneq V(F)$ of roots. Then the \emph{rooted $t$-blowup} of $F$ is the graph obtained by taking $t$ disjoint copies of $F$ and, for every $v\in R$, identifying the $t$ copies of $v$.

In this paper, we gave tight bounds for the extremal number of $K_{s,t}^k$. Notice that $K_{s,t}^k$ is the rooted $t$-blowup of $K_{s,1}^k$, where the roots are the leaves. This rooted graph is the $s$-legged spider with length vector $(k,\dots,k)$. It would be very interesting to extend the results to rooted blowups of spiders with other length vectors. To say more, we state the theorem of Bukh and Conlon we referred to in the introduction.

Let $F$ be a rooted graph with roots $R$. For a non-empty $S\s V(F)\setminus R$, we write $e_S$ for the number of edges incident to at least one vertex in $S$ and set $\rho_F(S)=\frac{e_S}{|S|}$. Moreover, we write $\rho(F)=\rho_F(V(F)\setminus R)$. We say that $F$ is \emph{balanced} if $\rho_F(S)\geq \rho(F)$ for every non-empty $S\s V(F)\setminus R$. Write $t\ast F$ for the rooted $t$-blowup of $F$. The result of Bukh and Conlon is as follows.

\begin{theorem}[Bukh--Conlon \cite{BC17}] \label{thmBC}
	Let $F$ be a balanced bipartite rooted graph with $\rho(F)>0$. Then there exists some $t_0=t_0(F)$ such that for all $t\geq t_0$, $$\ex(n,t\ast F)=\Omega(n^{2-\frac{1}{\rho(F)}}).$$
\end{theorem}

Now observe that the $s$-legged spider $F$ with length vector $(k,\dots,k)$ is balanced with $\rho(F)=\frac{sk}{s(k-1)+1}$, where the roots of $F$ are the leaves. Thus, Theorem \ref{thmBC} implies that $\ex(n,K_{s,t}^k)=\Omega(n^{1+\frac{s-1}{sk}})$ provided that $t$ is sufficiently large compared to $s$ and $k$, so Corollary~\ref{exactbip} follows.

More generally, it is not hard to see that an $s$-legged spider $S$ with length vector $(k_1,\dots,k_s)$ is balanced if and only if $k_1+\dots+k_s\geq (s-1)\max_{1\leq i\leq s}k_i$, where again the roots are the leaves. If this holds, then by Theorem \ref{thmBC}, $\ex(n,t\ast S)=\Omega(n^{1+\frac{s-1}{k_1+\dots+k_s}})$ for $t$ sufficiently large. This leads to the following natural conjecture.

\begin{conjecture} \label{spiderblowup}
	Let $s\geq 2$ and $1\leq k_1 \leq k_2\leq \dots \leq k_s$ be integers satisfying $k_1+\dots+k_s\geq (s-1)k_s$. Let $S$ be the rooted graph which is a spider with length vector $(k_1,\dots,k_s)$ and whose roots are the leaves. Then for any integer $t\geq 1$, $$\ex(n,t\ast S)=O(n^{1+\frac{s-1}{k_1+\dots+k_s}}).$$
\end{conjecture}

Our Theorem \ref{bipsub} proves this conjecture for $(k_1,\dots,k_s)=(k,\dots,k)$. The conjecture also holds for $k_1=1$, $k_2=\dots=k_s=k$ by Theorem 1.12 from \cite{CJL19}.

Some of the techniques in the present paper may be used to attack Conjecture \ref{spiderblowup}. More precisely, we have the following result.

\begin{lemma}
	Let $1\leq \ell_i\leq k_i$ be integers for each $i$. Assume that $\ell_i=1$ holds for at most one value of $i$. Let $G$ be a $K$-almost-regular graph on $n$ vertices with minimum degree $\d=\omega(1)$. Assume that $L$ is sufficiently large compared to $s$, $t$, $k_1,\dots,k_s$ and $K$ and that there are at least $\frac{n\d^{\ell_1+\dots+\ell_s}}{L}$ $L$-admissible but not $L$-good spiders with length vector $(\ell_1,\dots,\ell_s)$. Write $S$ for the spider with length vector $(k_1,\dots,k_s)$ and view it as a rooted graph with the roots being the leaves. Then $G$ contains $t\ast S$ as a subgraph.
\end{lemma}

\begin{proof2}
	Choose a set $\cS$ with the properties described in Lemma \ref{regularset}. Define $v_1,\dots,v_s$ as in the statement of Lemma \ref{connectthepaths}.
	We may repeatedly apply Lemma \ref{connectthepaths} to find $s$-legged spiders $S_1,\dots,S_t$, each with length vector $(k_1,\dots,k_s)$ and leaf vector $(v_1,\dots,v_s)$ such that $V(S_j)$ is disjoint from $(\bigcup_{1\leq i\leq j-1} V(S_i))\setminus \{v_1,\dots,v_s\}$. Then the union of these spiders is a copy of $t\ast S$.
\end{proof2}

\medskip

However, our other crucial ingredient, Lemma \ref{kpathsgood}, has not been generalised to the case where the forbidden subgraph is an arbitrary blowup of a spider (although, in \cite{CJL19} it is shown that one can use the blowup of the spider with length vector $(1,k,k\dots,k)$ in place of $K_{s,t}^k$).

\bibliographystyle{abbrv}
\bibliography{bibbipsub}
\bibliographystyle{plain}

\end{document}